\documentclass[10pt, letterpaper]{article}
\usepackage[utf8]{inputenc}
\usepackage{amsfonts, amsmath, amsthm, amssymb}
\usepackage{algpseudocode}
\usepackage{algorithm}
\usepackage{bbold}
\usepackage{graphicx}
\graphicspath{ {./img/} }
\usepackage{caption}
\usepackage{subcaption}
\usepackage{tabularx}
\usepackage{capt-of}

\usepackage[numbers]{natbib}

\algnewcommand{\IIf}[1]{\State\algorithmicif\ #1\ \algorithmicthen}
\algnewcommand{\EndIIf}{\unskip\ \algorithmicend\ \algorithmicif}

\usepackage{geometry}
\normalsize

\title{An Efficient Algorithm for Non-Negative Matrix Factorization with Random 
Projections}

\author{Gabriele Torre and Michael Graber\\
  \small{University of Applied Science and Art Northwestern Switzerland  (FHNW)}\\
  \small{Contact: firstname.lastname@fhnw.ch}\\
      	}

\begin{document}
      	\maketitle

\begin{abstract}
Non-negative matrix factorization (NMF) is one of the most popular
decomposition techniques for multivariate data. NMF is a core method for
many machine-learning related computational problems, such as data
compression, feature extraction, word embedding, recommender systems etc. 
In practice, however, its application is challenging for large datasets.
The efficiency of NMF is constrained by long data loading times, by large 
memory requirements and by limited parallelization capabilities.
Here we present a novel and efficient compressed NMF algorithm. Our
algorithm applies a random compression scheme to drastically reduce the
dimensionality of the problem, preserving well the pairwise distances 
between data points and inherently limiting the memory and communication load.
Our algorithm supersedes existing methods in speed. Nonetheless, it matches 
the best non-compressed algorithms in reconstruction precision.
\end{abstract}  

\section{Introduction}
Matrix factorizations constitute a fundamental pillar of numerous machine
learning methods. Essentially, matrix factorizations allow the approximate
decomposition of a multi-dimensional dataset into a linear combination of a 
limited set of components. 
The set of machine learning problems that can be tackled by the
use of matrix factorization are quite different in nature and span from
dimensionality reduction over blind source separation to the prediction of
ratings in collaborative filtering.\\
Principal Component Analysis (PCA) can be viewed as the most prominent matrix
factorization method. It allows to find mutually orthogonal, i.e. uncorrelated,
basis vectors along whose directions a given dataset shows decreasing variance.
By approximating the data points as linear combinations of a limited number 
of such orthonormal basis vectors, PCA allows to capture as much variance of 
the data as possible with a limited number of components. In this way, PCA can 
serve as an efficient \emph{data compression} technique.\\
Independent Component Analysis (ICA) allows to find basis vectors that are 
statistically independent. This represents a stronger constraint on the basis 
vectors than being uncorrelated. However, for certain applications 
this is desirable: ICA is well suited for the task of \emph{blind source
separation}, where different signal sources can be assumed to be statistically
independent \cite{hyvarinen2000independent}.\\
According to the conditions of orthogonality and independence imposed by PCA 
and ICA respectively, both of them find components that are mutually 
constrained. The support of the resulting component variables, however, is 
unconstrained, and often lead to basis vectors that cannot be easily 
interpreted. A natural property of many datasets is to have variables with 
non-negative support. For example, count-based measurements constitute a 
wide-ranging set of examples for non-negative datasets, including images based 
on photon counts, text document representations based on word counts or 
customer-product associations based on click or purchase counts.
Moreover, radiation spectra or network distance measurements represent
additional prominent non-negative dataset types.\\
Non-negative Matrix Factorization (NMF) allows for the decomposition of
non-negative datasets into two non-negative matrix factors.
Unlike PCA or ICA, NMF does not find components that are mutually orthogonal or
independent. It computes an approximation of the original dataset in terms of
additive linear combination of intrinsic non-negative features. This often
facilitates the discovery of naturally interpretable or even physically
meaningful data components where no further mutual geometric restrictions are 
necessary.\\
For this reason NMF has attracted interest during the last decade not only
in the fields of machine learning and data mining but also in application
domains, e.g. in Astronomy \cite{berne2010non}, Cosmology
\cite{zhang2016non}, Neuroscience \cite{maruyama2014detecting}, Recommender
Systems \cite{zhang2006learning}, and Document Clustering
\cite{shahnaz2006document}.\\
Despite its popularity, the application of standard NMF algorithms on large, 
high-dimensional data is restricted by their memory requirements. 
Standard NMF algorithms typically store the entire dataset in memory throughout 
the whole computational process. For this reason, the application of NMF 
methods on large datasets is often impractical.\\
To address these issues \cite{yu2014parallel} and \cite{kannan2016high}
recently introduced new NMF methods based on distributed and parallel
computation concepts. Furthermore, \cite{wang2010efficient} and
\cite{tepper2016compressed} showed that by using Semi-NMF methods 
\cite{ding2010convex} in combination with \emph{Random Projection} 
\cite{vempala2005random} the dimensionality of the problem can be drastically 
reduced and the underlying tasks solved much more efficiently.
It was proven that this approach is computationally more efficient and requires 
less memory when being compared to standard NMF algorithms 
\cite{lee1999learning}. 
Moreover, \emph{Random Projections} only marginally affect the quality of the 
NMF final results \cite{tepper2016compressed}.\\
In this paper we propose to incorporate a data compression scheme based on 
\emph{Random Projection} in the family of Hierarchical Alternating Least Square 
(HALS) NMF methods. HALS and FastHALS \cite{cichocki2009fast, 
cichocki2007hierarchical} are currently the fastest existing NMF algorithms. 
We will give a theoretical derivation of the methods and provide an 
empirical evaluation of their accuracy and numerical performances.\\
The structure of the paper is the following: In Section 2 we provide an 
introduction to the decomposition problem for non-negative datasets together
with some of the most popular NMF methods. In Section 3 we introduce the Random 
Projection technique as a data compression method. The detailed derivation of 
our algorithm is provided in Section 4. In Section 5 we present a set of 
experimental results obtained by means of the new NMF method and Section 
6 is devoted to our final discussions and considerations about the provided 
results.

\section{Model and Problem Statements}
For a given non-negative matrix $\mathbf{X} \in \mathbb{R}_+^{d\times n}$
composed of $d$ datapoints of dimensionality $n$, NMF seeks to 
identify the two non-negative factors 
$\mathbf{A} \in \mathbb{R}_+^{d \times k}$ and 
$\mathbf{B} \in \mathbb{R}_+^{n\times k}$, that provide a low-rank 
approximation of the form:
\begin{eqnarray}\label{eq.2.1}
  \mathbf{X} \sim \mathbf{A} \mathbf{B}^{T}~,
\end{eqnarray}
where $k \ll min(d,n)$ represents the desired number of components.
We introduce the cost function $J(\mathbf{A},\mathbf{B}|\mathbf{X})$ as the 
squared Euclidean distance (\emph{Frobenius norm}) between the data 
$\mathbf{X}$ and the computed approximation $\mathbf{A} \mathbf{B}^{T}$.
Starting from $J(\mathbf{A},\mathbf{B}|\mathbf{X})$, the global optimization 
problem of finding $\mathbf{A}$ and $\mathbf{B}$ given $\mathbf{X}$ can be 
expressed as:
\begin{equation}\label{eq2.2}
  \underset{\mathbf{A}, \mathbf{A}\geq0; \mathbf{B}, \mathbf{B}\geq0}
    {\text{argmin}} 
    J(\mathbf{A}, \mathbf{B} ~|~ \mathbf{X}) = \frac{1}{2}
        \| \mathbf{X} - \mathbf{A} \mathbf{B}^{T} \|_F^2~.
\end{equation}
Considering both $\mathbf{A}$ and $\mathbf{B}$ as variables of equation 
(\ref{eq2.2}), it has been proven that $J(\mathbf{A},\mathbf{B}|\mathbf{X})$ 
is a non-convex function \cite{lee2001algorithms}.
However, we can find two convex sub-problems if we consider $\mathbf{A}$ and 
$\mathbf{B}$ individually. 
Accordingly, a \emph{block-coordinate descent} approach 
\cite{lee2001algorithms} allows to compute values for $\mathbf{A}$ and 
$\mathbf{B}$ that correspond to a \emph{local} minimum of 
$J(\mathbf{A},\mathbf{B} ~|~\mathbf{X})$.

\subsection{Block Coordinate Descent Methods}\label{sec.2.1}
Generally, the iterative scheme adopted by block-coordinate descent algorithms
is to cyclically update \emph{blocks} of variables only, while keeping the 
remaining variables fixed. Assuming that the global constraints of the 
optimization problem can be decomposed into the Cartesian product of the 
constraint sets of each block of variables, the resulting sequence of 
block-coordinate updates is guaranteed to converge to a stationary point
\cite{bertsekas1995dynamic}.

In the case of NMF, the overall non-negativity constraints are in fact the
Cartesian product of the non-negativity constraints on the individual
variables. Hence, the NMF optimization problem can be tackled in a
block-coordinate descent approach. 
The most elemental block coordinate descent approach to NMF uses $\mathbf{A}$ 
and $\mathbf{B}$ as coordinate block. The resulting optimization technique is 
sketched in Algorithm 1.

\begin{algorithm}\label{A1}
\caption{NMF block-coordinate descent approach.}
\begin{algorithmic}[1]
\For{iterations ($i$)}
\State \small{$\mathbf{A}^{(i+1)} \leftarrow 
\underset{\mathbf{A},\mathbf{A}\geq0}{\text{argmin}} ~ 
    J(\mathbf{A}^{(i)},\mathbf{B}^{(i)} ~|~\mathbf{X}) = 
    \frac{1}{2}\| \mathbf{X} - \mathbf{A}^{(i)} \mathbf{B}^{(i)T} \|_F^2$}
    \Comment{$\mathbf{B}$ fixed}
\State \small{$\mathbf{B}^{(i+i)} \leftarrow 
\underset{\mathbf{B},\mathbf{B}\geq0}{\text{argmin}} ~ 
    J(\mathbf{A}^{(i+1)},\mathbf{B}^{(i)}~|~\mathbf{X}) = 
    \frac{1}{2}\| \mathbf{X} - \mathbf{A}^{(i+1)} \mathbf{B}^{(i)T} \|_F^2 $}
    \Comment{$\mathbf{A}$ fixed}
    \IIf{$\text{convergence criterion is reached}$} $\text{Stop iterations}$\EndIIf
\EndFor
\end{algorithmic}
\end{algorithm}

NMF methods which adopt this optimization technique are, e.g., the 
\emph{Multiplicative Updates rule} \cite{lee1999learning}, the 
\emph{Active-Set-Like} method \cite{kim2008nonnegative} or 
\emph{Projected Gradient Descent} NMF \cite{lin2007projected}.
However, these methods are \textit{computationally expensive} and associated 
to a \textit{slow convergence rate} \cite{kim2014algorithms}. 
In fact, they are characterized by large memory needs that scale with 
$\mathcal{O}(dn+dk+nk)$ and involve a large number of Floating Point 
Operations (FLOPS) per iterative update ($\sim \mathcal{O}(2dnk)$).

\subsection{Hierarchical Optimization Methods}\label{sec.2.2}
The Hierarchical Alternating Least Squares (HALS) method for NMF was 
originally proposed by \cite{cichocki2007hierarchical, cichocki2009fast} as an 
improvement of the Alternating Least Squares (ALS) method 
\cite{cichocki2007regularized}.
It consists of a block-coordinate descent method with single component vectors 
as coordinate blocks. In this context, the cost function from equation 
(\ref{eq2.2}) can be modified to define a set of cost functions:
\begin{equation}\label{cost2}
J(\mathbf{a}_j, \mathbf{b}_j ~|~\mathbf{X}_j) = 
\| \mathbf{X}_j - \mathbf{a}_j\mathbf{b}_j^T \|_F^2 \qquad \forall j \in [1,\cdots,k]~,
\end{equation}
one for each component vector $\mathbf{a}_j$ and $\mathbf{b}_j$ such that
$\mathbf{A} = [\mathbf{a}_1,\mathbf{a}_2,\cdots,\mathbf{a}_k]$ and
$\mathbf{B} = [\mathbf{b}_1,\mathbf{b}_2,\cdots,\mathbf{b}_k]$ and where
$\mathbf{X}_j$ is: 
\begin{equation}\label{eq.2.2_3}
  \mathbf{X}_{j} = \mathbf{X} - \mathbf{A}\mathbf{B}^T + \mathbf{a}_j
  \mathbf{b}_j^T \qquad 
  \forall j \in [1,\cdots,k]~,
\end{equation}
where the computed expectation of $\mathbf{X}$, provided by the $j^{th}$ 
components of $\mathbf{A}$ and $\mathbf{B}$ is added to the residual matrix.
The global optimization problem based on equation (\ref{cost2}) can be 
addressed by the set of iterative updates presented in Algorithm 2.

\begin{algorithm}
\caption{NMF Hierarchical Alternating Least Squares approach.}
\label{alg:A2}
\begin{algorithmic}[1]
  \For{iterations ($i$)}
    \For{components ($j$)}
    \State \small{$\mathbf{a}_j^{(i+1)} \leftarrow
    \underset{\mathbf{a}_j,\mathbf{a}_j\geq0}{\text{argmin}} ~
    J(\mathbf{a}_j^{(i)},\mathbf{b}_j^{(i)} ~|~\mathbf{X}_j) = \frac{1}{2}
    \|\mathbf{X}_{j}-\mathbf{a}_j^{(i)} \mathbf{b}_j^{(i)T} \|_F^2$}
        \Comment{$\mathbf{b}_j$ fixed}
    \State \small{$\mathbf{b}_j^{(i+1)} \leftarrow
    \underset{\mathbf{b}_j,\mathbf{b}_j\geq0}{\text{argmin}} ~
    J(\mathbf{a}_j^{(i+1)},\mathbf{b}_j^{(i)} ~|~\mathbf{X}_j) = \frac{1}{2}
    \|\mathbf{X}_{j}-\mathbf{a}_j^{(i+1)} \mathbf{b}^{(i)T}_j \|_F^2$}
        \Comment{$\mathbf{a}_j$ fixed}
    \EndFor
    \IIf{$\text{convergence criterion is reached}$} $\text{Stop iterations}$\EndIIf
    \EndFor
\end{algorithmic}
\end{algorithm}
The HALS method provides alternating updates on a single component level 
which typically result in higher convergence rate and better data 
approximation compared to the methods introduced on Section \ref{sec.2.1}
\cite{kim2014algorithms}.
On the other hand, the iterative update rules of HALS are computationally 
more expensive if compared to the block coordinate descent based method, 
with a number of FLOPS per iteration $\sim \mathcal{O}(8dnk)$.  
This computational limit was partially overcome by the introduction of 
FastHALS \cite{cichocki2007regularized}, a more efficient update rule for 
HALS which scales as $\mathcal{O}(4dnk)$.
However, a critical issue remains: For large datasets the HALS algorithm 
still requires the entire data matrix $\mathbf{X}$ to be held in memory in 
order to be fast. Hence the memory consumption still scales with 
$\mathcal{O}(dn+dk+nk)$.

\section{Random Projection}\label{sec.3}
Random projection is a dimensionality reduction technique for datapoints 
lying in a Euclidean space. It is commonly adopted to reduce the problem of 
managing and manipulating large datasets for techniques such as PCA, 
Singular Value Decomposition, Manifold Learning \cite{wang2010efficient} 
and NMF \cite{tepper2016compressed, bingham2001random}. 
Its significance is mostly due to the \emph{Johnson-Lindenstrauss lemma} which
proves that random projections well preserve the pairwise Euclidean distances 
between datapoints \cite{dasgupta2003elementary}.
For this exact same reason, NMF methods with random projection provide data 
approximation results almost as accurate as their standard uncompressed 
counterparts \cite{tepper2016compressed}.
In this section we provide a general introduction to the random projections
technique while the application of this method to the context of NMF will be 
discussed in section \ref{sec.3.1}.\\

The random projections are characterized by a simple computational scheme. 
Let us define $r \in \mathbb{N}$ as the rank of the given data matrix
$\mathbf{X} \in \mathbb{R}^{d \times n}$, with $r_{ov} \in \mathbb{N}$ as an
oversampling parameter and $\mathbf{\Omega} \in \mathcal{N}(0,1)^{n \times
(r+r_{ov})}$ as a Gaussian random matrix whose entries are standard normal
random variables. 
The structured random projection operation \cite{tepper2016compressed} can 
be defined as:
\begin{equation}\label{eq.3.4}
  \hat{\mathbf{X}} = \mathbf{Q}^T\mathbf{X}
\end{equation}
where $\mathbf{Q}:= \mathbb{R}^{d \times n} \rightarrow 
\mathbb{R}^{(r+r_{ov}) \times n}$ is the random projection operator that maps 
the datapoints from the original dataspace onto the randomly selected subspace.
The orthogonal matrix $\mathbf{Q}$ is defined by the QR decomposition of the 
data transformation $P(\mathbf{X})$:
\begin{equation}\label{eq.3.5}
  P(\mathbf{X}) = (\mathbf{X}\mathbf{X}^T)^w \mathbf{X} \mathbf{\Omega},
\end{equation}
where $w \in \mathbb{N}$ is the power iterations parameter 
\cite{halko2011finding}.
For large input matrices $\mathbf{X}$, the singular vectors associated with
small singular values will interfere with the calculation of the random
projection matrix. Thanks to the factor $(\mathbf{X}\mathbf{X}^T)^w$ the decay
rate of the singular power spectrum $\sigma_j (\mathbf{X})$ will be increased
to:
\begin{equation}
    \sigma_j \left( (\mathbf{X}\mathbf{X}^T)^w \mathbf{X} \right) = 
    \sigma_j \left(\mathbf{X} \right)^{2w+1}~\text{with}~j = 1,2,3,... ~~.
\end{equation}
This preserves the singular vectors of $\mathbf{X}$ while rendering the larger
singular values more dominant for the definition of $\mathbf{Q}$.
In the following section we will introduce the state of the art techniques 
which incorporate random projections into the block coordinate descent approach 
for NMF.

\subsection{Random Projection and Semi-NMF}\label{sec.3.1}
According to the definition of the operator $\mathbf{Q}$ given in equation 
(\ref{eq.3.4}), random projection operators do not preserve the property of 
non-negativity of a non-negative dataset $\mathbf{X}$.
For this reason, NMF techniques with random projection are addressed in terms 
of two different objective functions which are defined in terms of the two non 
negative factor matrices $\mathbf{A}$ and $\mathbf{B}$.\\
Let us define the two random projection matrices 
$\mathbf{L} \in \mathbb{R}^{d \times (r+r_{ov})}$ and 
$\mathbf{R} \in \mathbb{R}^{(r+r_{ov}) \times n}$ respectively 
from $P(\mathbf{X})$ and $P(\mathbf{X}^T)$ and assuming that 
$min(d,n) \gg (r+r_{ov}) > k$. The NMF problem with random projection can be 
defined as the alternating optimization problem, as shown in Algorithm 3, 
where the optimization problems of line 2 and 3 are addressed in terms of 
the Semi-NMF method \cite{ding2010convex}.

\begin{algorithm}
  \begin{algorithmic}[1]
    \caption{semi-NMF with random projection}\label{alg:A3}
      \For{iterations ($i$)}
      \State \small{$ \mathbf{A}^{(i+1)} \leftarrow 
      \underset{\mathbf{A},\mathbf{A}{\geq0}}{\text{argmin}} ~
    J(\mathbf{A}^{(i)},\mathbf{B}^{(i)}\mathbf{R}^T | \mathbf{X}\mathbf{R}^T)= 
      \frac{1}{2}\| \mathbf{X}\mathbf{R}^T - 
      \mathbf{A}^{(i)}\mathbf{B}^{(i)T}\mathbf{R}^T \|_F^2$}
      \Comment{$\mathbf{B}$ fixed}
      \State \small{$ \mathbf{B}^{(i+1)} \leftarrow 
      \underset{\mathbf{B},\mathbf{B}{\geq0}}{\text{argmin}} ~
    J(\mathbf{L}^T\mathbf{A}^{(i+1)},\mathbf{B}^{(i)}|\mathbf{L}^T\mathbf{X})=
      \frac{1}{2}\| \mathbf{L}^T\mathbf{X} - 
      \mathbf{L}^T\mathbf{A}^{(i+1)}\mathbf{B}^{(i)T} \|_F^2$}
      \Comment{$\mathbf{A}$ fixed}
      \IIf{$\text{convergence criterion is reached}$} $\text{Stop iterations}$\EndIIf
    \EndFor
  \end{algorithmic}
\end{algorithm}

\noindent  According to Algorithm \ref{alg:A3}, only the currently updated 
and not projected factor matrix is assumed to be non-negative.

\section{HALS with Random Projection}\label{sect4}
In this section we formally introduce two new NMF block coordinate descent 
approaches incorporating random projections, based on the HALS and FastHALS
methods. According to the definition of $\mathbf{L}$ and $\mathbf{R}$, given
in section \ref{sec.3.1}, Algorithm \ref{alg:A4} shows how random projections
can be included in the computational scheme described by Algorithm
\ref{alg:A2}, where $\hat{\mathbf{X}} = \mathbf{L}^T\mathbf{X}$ and 
$\check{\mathbf{X}} = \mathbf{X}\mathbf{R}^T$ are the \emph{left-}
and \emph{right-projected} datasets respectively, 
$\hat{\mathbf{A}}=\mathbf{L}^T\mathbf{A}$ and $\check{\mathbf{B}}=
\mathbf{R}\mathbf{B}$ are the compressed factors with components 
$\hat{\mathbf{a}}_j = \mathbf{L}^T \mathbf{a}_j$ and 
$\check{\mathbf{b}}_j = \mathbf{R}\mathbf{b}_j$ and 
$\hat{\mathbf{X}}_{j}$ and $\check{\mathbf{X}}_{j}$ are defined as:
\begin{equation}\begin{split}\label{eq.4_6}
 \hat{\mathbf{X}}_{j}& = \hat{\mathbf{X}} - \hat{\mathbf{A}}\mathbf{B}^T + 
  \hat{\mathbf{a}}_j \mathbf{b}^T_j\\
 \check{\mathbf{X}}_{j}& = \check{\mathbf{X}} - \mathbf{A}\check{\mathbf{B}}^T 
  + \mathbf{a}_j\check{\mathbf{b}}_j^T.\\
\end{split}\end{equation}

\begin{algorithm}
\begin{algorithmic}[1]
\caption{Hierarchical Alternating Least Squares with Random Projection}
\label{alg:A4}
\For{iteration ($i$)}
\For{components ($j$)}
\State \small{$\mathbf{a}_j^{(i+1)} \leftarrow 
  \underset{\mathbf{a}_j,\mathbf{a}_j\geq0}{\text{argmin}} ~ 
    J(\mathbf{a}_j^{(i)},\check{\mathbf{b}}_j^{(i)} ~|~\check{\mathbf{X}}_j) = 
    \frac{1}{2}
    \|\check{\mathbf{X}}_j-\mathbf{a}^{(i)}_j\check{\mathbf{b}}^{(i)T}_j\|_F^2$
    } \label{A4_L3} 
\Comment{$\mathbf{B}$ fixed}
  
\State \small{$\mathbf{b}_j^{(i+1)} \leftarrow 
  \underset{\mathbf{b}_j,\mathbf{b}_j\geq 0}{\text{argmin}} ~ 
  J(\hat{\mathbf{a}}_j^{(i+1)}\mathbf{b}_j^{(i)} ~|~\hat{\mathbf{X}}_j) = 
  \frac{1}{2}
  \|\hat{\mathbf{X}}_j-\hat{\mathbf{a}}^{(i+1)}_j\mathbf{b}_j^{(i)T}\|_F^2 $}
  \label{A4_L4}
\Comment{$\mathbf{A}$ fixed}
\EndFor
\IIf{$\text{convergence criterion is reached}$} $\text{Stop iterations}$\EndIIf
\EndFor
\end{algorithmic}
\end{algorithm}

The explicit form of the iterative update rules in line 3 and 4 of Algorithm 
\ref{alg:A4} can be found by computing the local gradient of the two cost 
functions $J(\mathbf{a}_j,\check{\mathbf{b}}_j ~|~\check{\mathbf{X}}_{j})$ and 
$J(\hat{\mathbf{a}}_j,\mathbf{b}_j ~|~\hat{\mathbf{X}}_{j})$ with respect to 
the unknown factor vectors $\mathbf{a}_j$ and $\mathbf{b}_j$ as:
\begin{equation}\label{hals_a}
  \frac{\partial J(\mathbf{a}_j,\check{\mathbf{b}}_j ~|~\check{\mathbf{X}}_{j})}
  {\partial \mathbf{a}_j} =
  (\mathbf{a}_j\check{\mathbf{b}}_j^T\check{\mathbf{b}}_j - 
  \check{\mathbf{X}}_j \check{\mathbf{b}}_j),
\end{equation}
\begin{equation}\label{hals_b}
  \frac{\partial J(\hat{\mathbf{a}}_j,\mathbf{b}_j ~|~\hat{\mathbf{X}}_{j})}
  {\partial \mathbf{b}_j} =
  (\mathbf{b}_j\hat{\mathbf{a}}^T_j\hat{\mathbf{a}_j} - 
  \hat{\mathbf{X}}_{j}^T \hat{\mathbf{a}}_j).
\end{equation}
By equating the gradient components to zero, the iterative update scheme 
for the HALS-RP method follows directly from equation \eqref{hals_a} and 
\eqref{hals_b} as:
\begin{equation}\label{hals_update}
  \mathbf{a}_j \leftarrow 
  \left[ \frac{\check{\mathbf{X}}_{j} \check{\mathbf{b}}_j}
    {\check{\mathbf{b}}_j^T\check{\mathbf{b}}_j} \right]_+,~
  \mathbf{b}_j \leftarrow 
  \left[ \frac{\hat{\mathbf{X}}_{j}^T \hat{\mathbf{a}}_j}
  {\hat{\mathbf{a}}_j^T\hat{\mathbf{a}}_j}\right]_+~;
\end{equation}
where the non-negativity of $\mathbf{a}_j$ and 
$\mathbf{b}_j$ is imposed after each iterative update by setting the negative
values of the solution to zero.\\

According to the derivation of the FastHALS method presented in 
\cite{cichocki2007hierarchical}, we introduce an alternative and computationally 
more efficient update rule for the HALS-RP algorithm, named FastHALS-RP. 
This new method is defined by including the explicit expressions for 
$\hat{\mathbf{X}}$ and $\check{\mathbf{X}}$ from equation (\ref{eq.4_6}) in 
equations (\ref{hals_update}) as:

\begin{equation}\label{fastHALS_a}
\begin{split}
\mathbf{a}_j \leftarrow & 
  \left[ \mathbf{a}_j + 
    \frac{ \check{\mathbf{X}}\check{\mathbf{b}}_j - 
        \mathbf{A}\check{\mathbf{B}}^T\check{\mathbf{b}}_j }
    { \check{\mathbf{b}}_j^T\check{\mathbf{b}}_j } 
        \right]_+ = 
  \left[ \mathbf{a}_j + 
    \frac{ \left(\check{\mathbf{X}}\check{\mathbf{B}}\right)_{j} - 
    \mathbf{A}\left(\check{\mathbf{B}}^T\check{\mathbf{B}}\right)_{j} }
    {\left(\check{\mathbf{B}}^T\check{\mathbf{B}}\right)_{jj}}
    \right]_+,
    \\ 
\mathbf{b}_j \leftarrow &
  \left[ \mathbf{b}_j + 
    \frac{ \hat{\mathbf{X}}^T\hat{\mathbf{a}}_j - 
    \mathbf{B}\hat{\mathbf{A}}^T\hat{\mathbf{a}}_j }
          { \hat{\mathbf{a}}_j^T\hat{\mathbf{a}}_j }    
  \right]_+ = 
  \left[ \mathbf{b}_j + 
    \frac{ \left(\hat{\mathbf{X}}^T\hat{\mathbf{A}}\right)_{j} - 
    \mathbf{B}\left(\hat{\mathbf{A}}^T\hat{\mathbf{A}}\right)_{j}}
    { \left(\hat{\mathbf{A}}^T\hat{\mathbf{A}}\right)_{jj}}\right]_+~;
\end{split}
\end{equation}
where the normalization of the factor vectors is iteratively imposed with
$\mathbf{a}_{j} = \mathbf{a}_{j} / \|\mathbf{a}_j\|_2$.

\subsection{Sparsity and Smoothness constraints}
In order to enforce the properties of smoothness and sparsity of the
factor matrix $\mathbf{B}$, we introduce two additional penalty terms to 
the likelihood function presented in line 4 of Algorithm \ref{alg:A4} as: 
\begin{equation}\label{obj_fun_sp_sm}
  J(\hat{\mathbf{a}}_j,\mathbf{b}_j ~|~\hat{\mathbf{X}}_{j}) = 
  \frac{1}{2} \|\hat{\mathbf{X}}_{j} - 
  \hat{\mathbf{a}}_j \mathbf{b}^T_j \|_F^2 + 
  \alpha \|\mathbf{b}_j\|_1 +
    \frac{\beta}{2} \|\mathbf{b}_j\|_2^2,
\end{equation}
where $\alpha \in \mathbb{R}_+$ and $\beta \in \mathbb{R}_+$ are 
the two parameters regulating the sparsity and smoothness levels of 
$\mathbf{b}_j$ while $\| \cdot \|_1$ is the $L_1$ norm.
Following the same derivation scheme presented in section \ref{sect4}, 
the HALS-RP update rules for $\mathbf{a}_j$ and $\mathbf{b}_j$ are:
\begin{equation}\label{hals_constr}
\begin{split}
  \mathbf{a}_j \leftarrow &
  \left[ \frac{\check{\mathbf{X}}_{j} \check{\mathbf{b}}_j}
  {\check{\mathbf{b}}_j^T\check{\mathbf{b}}_j} \right]_+,\\
  \mathbf{b}_j \leftarrow & 
  \frac{\left[\hat{\mathbf{X}}_{j}^T \hat{\mathbf{a}}_j - 
    \alpha \mathbb{1}_n\right]_+ }
    {\hat{\mathbf{a}}_j^T\hat{\mathbf{a}}_j + \beta},
\end{split}
\end{equation}
where $\mathbb{1}_n$ denotes the vector made of all unitary entries with 
length $n$. Analogously, imposing the same constraints for the vector
$\mathbf{b}_j$ computed by the FastHALS-RP method, we get the following 
update rules:
\begin{equation}
\begin{split}
\mathbf{a}_j \leftarrow &
    \left[ \mathbf{a}_j + \frac{ \mathbf{H}_{j} -
        \mathbf{A}\mathbf{G}_{j} }
        {\mathbf{G}_{jj}}
     \right]_+,\\
    \mathbf{b}_j \leftarrow &
     \left[\frac{ \mathbf{b}_j \hat{\mathbf{a}_j}^T\hat{\mathbf{a}}_j + 
      \hat{\mathbf{X}}^T\hat{\mathbf{a}}_j -  
      \mathbf{B}\hat{\mathbf{A}}^T\hat{\mathbf{a}}_j^T  -
    \alpha\mathbb{1}_N}{\hat{\mathbf{a}}_j^T\hat{\mathbf{a}}_j + 
    \beta}\right]_+ =  \\
    & \left[ \frac{ \mathbf{b}_j \mathbf{W}_{jj} + \mathbf{P}_{j} - 
    \mathbf{B}\mathbf{W}_{j} - \alpha\mathbb{1}_N}
    {\mathbf{W}_{jj} + \beta}\right]_+,
  \end{split}
\end{equation}
where $\mathbf{H} = \check{\mathbf{X}}\check{\mathbf{B}}$, 
$\mathbf{G} = \check{\mathbf{B}}^T \check{\mathbf{B}}$, 
$\mathbf{W} = \hat{\mathbf{A}}^T\hat{\mathbf{A}}$ and 
$\mathbf{P} = \hat{\mathbf{X}}^T\hat{\mathbf{A}}$.

\subsection{Numerical Complexity and Memory Consumption}\label{sec.4.2}
This section is aimed at evaluating the theoretical performance expectations 
of the different NMF algorithms presented in the previous sections. 
In particular, we are interested in assessing the improvements provided by 
random projections in terms of numerical complexity and memory consumption. 
Table \ref{table.4.1} shows the performance comparison between the 
Multiplicative Update method (MU), the Hierarchical Least Squares method 
(HALS) and its more efficient implementation (FastHALS) with their 
variants using random projections i.e., MU-RP, HALS-RP and FastHALS-RP.

\begin{table}[h!]
\begin{small}
\begin{center}
\begin{tabular}{ | l | c | c |}
\hline
Method & Num. Complexity & Memory Consumption\\
\hline
\hline
MU & $\mathcal{O}(8dnk)$ & $\mathcal{O}(dn + dk + nk ) $ \\
\hline
HALS & $\mathcal{O}(8dnk)$ & $\mathcal{O}(dn + dk + nk) $\\
\hline
FastHALS & $\mathcal{O}(4dnk)$ & $\mathcal{O}(dn + dk + nk) $\\
\hline
\hline
MU-RP & $\mathcal{O}(4dk(r+r_{ov}))$ & $\mathcal{O}((2(r+r_{ov})+k)(d+n))$\\
\hline
HALS-RP & $\mathcal{O}(4dnk(r+r_{ov}))$ & $\mathcal{O}((2(r+r_{ov})+k)(d+n))$\\
\hline
FastHALS-RP & $\mathcal{O}(2dk(r+r_{ov}))$ & $\mathcal{O}((2(r+r_{ov})+k)(d+n))$\\
\hline
\end{tabular}
\end{center}
\end{small}
\caption{Numerical complexity and memory requirements for MU, HALS, FastHALS, 
MU-RP, HALS-RP and FastHALS NMF. All those values are computed under the 
assumption $d \gg (n,r,k)$.}
\label{table.4.1}
\end{table}

According to Table \ref{table.4.1}, random projections reduces the numerical 
complexity associated to MU and FastHALS by a factor $2(r+r_{ov}) / n$, while 
it increases by a factor $(r+r_{ov}) / 2$ for the HALS method. 
This discrepancy is a direct consequence of the methods' different computational 
arrangements.
In particular, HALS-RP requires $k$ different random projection steps within 
every iterative loop (one for each component $\mathbf{a}_j$ and $\mathbf{b}_j$),
leading to a strong increment of the numerical complexity in its iterative update 
step. In conclusion, random projections provide a reduction of the memory 
consumptions from $\mathcal{O}(dn+dk+nk)$ to $\mathcal{O}((2(r+r_{ov})+k)(d+n))$ 
for all the presented NMF methods.

\section{Experiments}
This section is aimed at providing an overview of the numerical and 
computational properties of the NMF methods based on two real-world 
applications: a factorization of the \emph{Olivetti faces dataset} and the
\emph{20 Newsgroups dataset} which represent dense and sparse data
respectively.
Our analysis is structured as follows: First, we investigate some of 
numerical and computational properties of the NMF methods, such as, 
convergence rate, reconstruction errors, numerical complexity and memory 
consumption. Second, we evaluate how the solutions computed by FastHALS and 
FastHALS-RP are influenced by different choices for the number 
of components $k$, the random projection parameter $w$ and the sparsity and 
smoothness parameters, respectively $\alpha$ and $\beta$. 
For the assessment of the performance in the second step, we discuss measurable
quantities, e.g. the final data approximation error and the sparsity level of
$\mathbf{B}$.
The measured quantities presented here were estimated as the median values 
over multiple independent runs of every algorithm to avoid biases in the 
results and to assess the stability of the methods against the initial random 
initialization of the factor matrices.

\subsection{Olivetti Faces dataset}\label{sec5.1}
\begin{figure}[htbp]
  \includegraphics[width=0.95\textwidth]{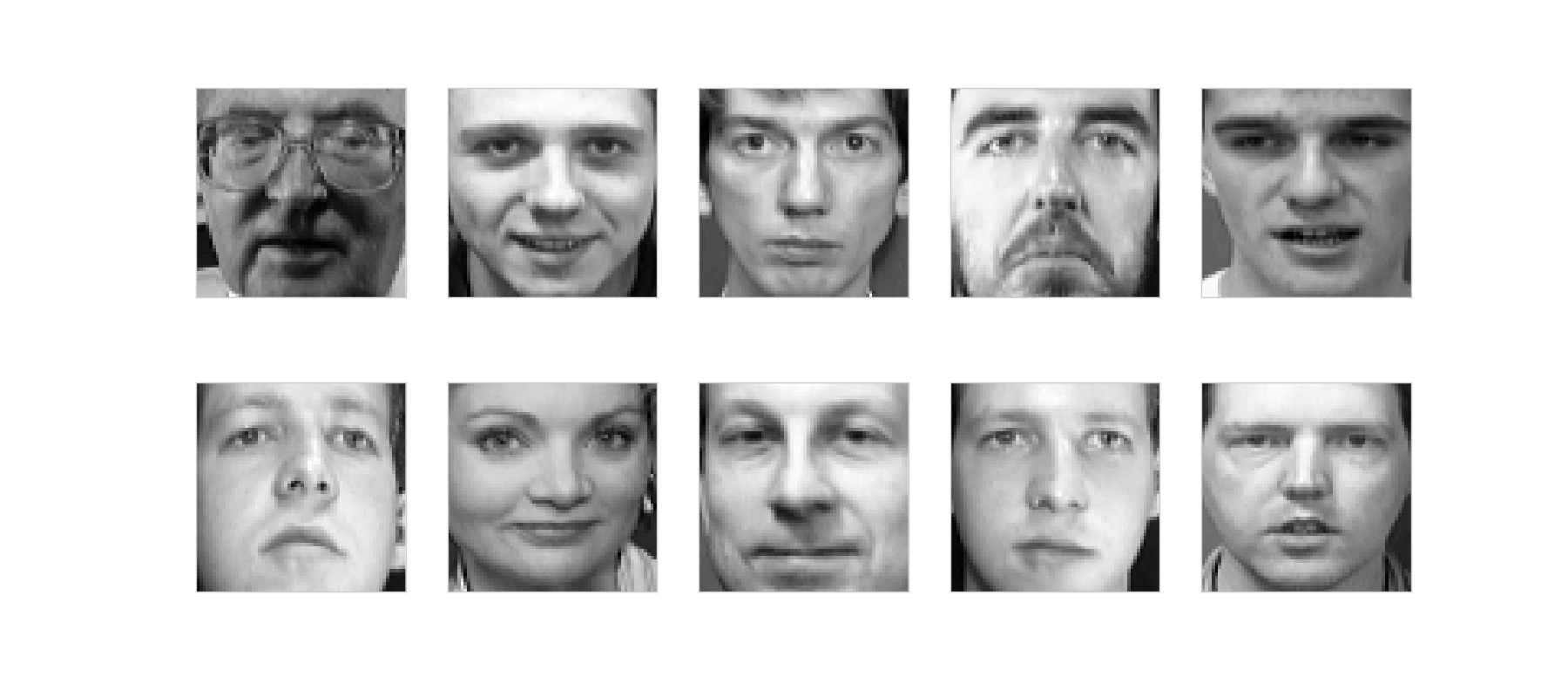}
  \caption{Ten samples from the Olivetti faces dataset of images.}
  \label{fig.5.1}
\end{figure}

The Olivetti faces dataset \cite{olivettioracle} is composed of 400 images
originally collected for the validation of face recognition algorithms
\cite{samaria1994face}. 
Images, $64\times64$ pixels, are quantized to 256 gray scale levels and show 
the faces of 40 different subjects under varying light conditions and with
different facial expressions (Figure \ref{fig.5.1}).

\begin{figure}[htbp]
      \includegraphics[width=0.95\textwidth]{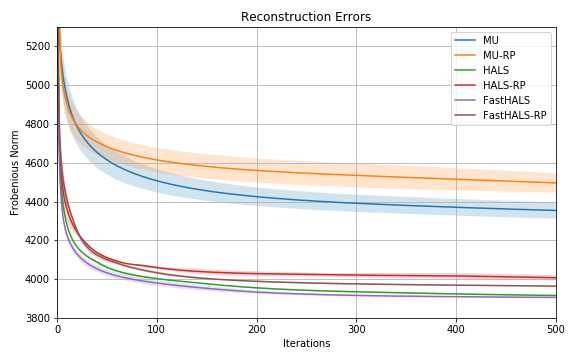}
      \caption{Evolution of the reconstruction errors over iterations for the
      MU, MU-RP, HALS, HALS-RP, FastHALS, FastHALS-RP NMF methods when applied
      to the Olivetti faces dataset.}
      \label{fig5.2}
\end{figure}

The first set of results presented in this section is based on runs over 500 
iterations of the MU, MU-RP, HALS, HALS-RP, FastHALS and FastHALS-RP NMF 
algorithms applied to the Olivetti faces dataset. We choose the number of 
components $k=20$ for all the NMF methods while we keep $r+r_{ov}=25$ and
$w = 4$ fixed for the MU-RP, HALS-RP and FastHALS-RP algorithms. Figure
\ref{fig5.2} shows the iterative evolution of the NMF data approximation 
error computed as $\frac{1}{2}\|\mathbf{X}-\mathbf{A}\mathbf{B}^T\|_F^2$. 
In comparison with the other NMF methods, MU and MU-RP performs with higher 
reconstruction errors and smaller iterative convergence rates which can be 
interpreted as a direct consequence of the different computational granularity 
between the MU and the HALS based algorithms. Moreover, by considering the 
plot in Figure \ref{fig5.2}, it is not clear if MU and MU-RP will asymptotically 
converge to the same data approximation precision level of the other NMF methods.
In table \ref{table5.2} we present the numerical complexity, memory 
consumption and average time per iteration values estimated for this 
application. According to these results, FastHALS-RP, together with MU-RP, 
outperform the other NMF methods in terms of average time per update. While 
FastHALS-RP has the lowest numerical complexity. In particular, FastHALS-RP 
allow for a numerical complexity level which is about the $40\%$ lower with 
respect to the FastHALS method.

As already discussed in section \ref{sec.4.2}, these results illustrate that
in the HALS case, random projections come with additional computational cost.
In fact, HALS-RP is computationally more expensive than HALS by more than one 
order of magnitude. Finally, the last column of table \ref{table5.2} shows 
that random projections provide a reduction of the NMF algorithms memory 
consumption of about the $80\%$.

\begin{table}[htbp]
\centering
\begin{tabular}{ |l||c|c||c| }
        \hline
        Method       & Complexity     & Time per update  & Memory          \\
                     & [$10^6$ flops] & [seconds]        & [$10^4$ floats] \\   
        \hline
        MU           & $270$          & 0.031            & $170$           \\
        \hline
        MU-RP        & $13$           & 0.023            & $31$            \\
        \hline
        HALS         & $260$          & 0.332            & $170$           \\
        \hline
        HALS-RP      & $3200$         & 0.575            & $31$            \\
        \hline
        FastHALS     & $14$           & 0.035            & $170$           \\
        \hline
        FastHALS-RP  & $8.6$          & 0.023            & $31$            \\
        \hline
    \end{tabular}
    \caption{NMF methods numerical complexity, memory consumption and 
            time per update estimated for the Olivetti faces applications.}
            \label{table5.2}
\end{table}

The second set of assessments included in this section is aimed at estimating 
how different choices for the parameters $k$, $w$, $\alpha$ and $\beta$ affect 
the approximation error of the FastHALS-RP method. Additionally, we investigate
the effect on sparsity and smoothness properties of the computed factor
$\mathbf{B}$.

\begin{figure}[htbp]
  \includegraphics[width=0.95\textwidth]{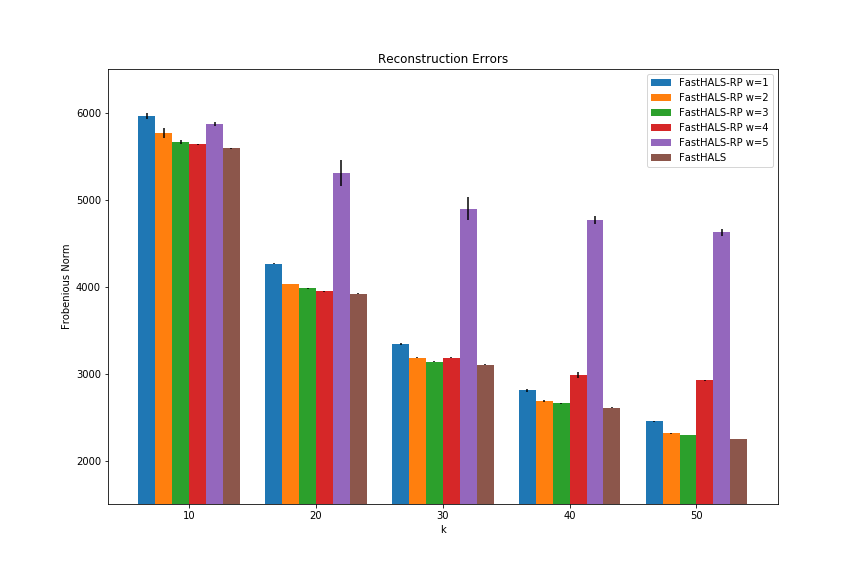}
  \caption{Comparison between the approximation errors provided by FastHALS 
           and FastHALS-RP for $k = [10,20,30,40,50]$ and $w = [1,2,3,4,5]$.}
  \label{fig5.3} 
\end{figure}

\begin{figure}[htbp]
\includegraphics[width=0.95\textwidth]{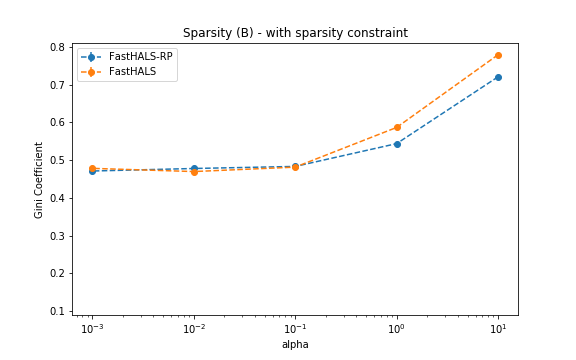}
\caption{Gini coefficients of the $\mathbf{B}$ matrix computed by FastHALS 
         and FastHALS-RP with sparsity parameter 
         $\alpha = [0.001,0.01,0.1,1,10]$, number of components $k=20$ 
         and random projection parameter $w = 3$.}\label{fig5.4a}
\end{figure}

\begin{figure}[htbp]
\includegraphics[width=0.95\textwidth]{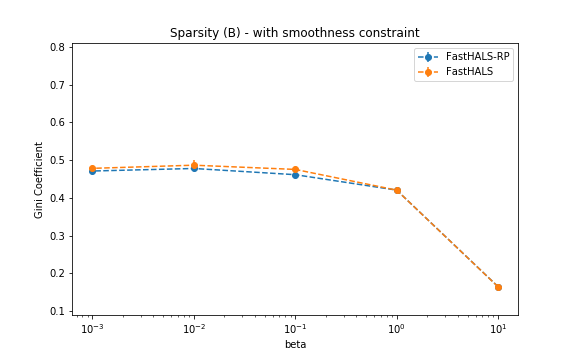}
\caption{Gini coefficients of the $\mathbf{B}$ matrix computed by FastHALS 
         and FastHALS-RP with smoothness parameter
         $\beta = [0.001,0.01,0.1,1,10]$, number of components $k=20$ 
         and random projection parameter $w = 3$.}\label{fig5.4b}
\end{figure}

Figure \ref{fig5.3} compares the data reconstruction quality provided by 
FastHALS-RP and FastHALS which appear to be strongly dependent on the choice 
made for the parameters $k$ and $w$. For the range of $k$ values spanned by 
this analysis, $w=3$ provides the most accurate results in terms of data
approximation precision.

Figures \ref{fig5.4a} and \ref{fig5.4b} show how the sparsity level of the
factor matrix $\mathbf{B}$ computed by FastHALS-RP and FastHALS is influenced 
by the parameters $\alpha$ and $\beta$. 
The sparsity level of $\mathbf{B}$ was estimated in terms of the Gini 
coefficient: 
\begin{equation}
   G = \frac{\sum_{i=1}^{n\cdot k} (2i - (n \cdot k) - 1) 
   \mathbf{B}^{S}_{i}}{(n\cdot k)\sum_{i=1}^{n\cdot k} \mathbf{B}^{S}_{i}},
\end{equation}
where $\mathbf{B}^{S}$ is the vector containing the values of the matrix
$\mathbf{B}$ sorted in ascending order. 
Figure \ref{fig5.4a} shows that the sparsity level of the factor matrix 
$\mathbf{B}$ grows consistently with the parameter $\alpha$. Moreover,
for $\alpha < 0.1$ the two methods provide similar solution in terms of 
the Gini Coefficients profiles, while they start to deviate from each 
other for $\alpha >= 0.1$ where FastHALS provides matrices $\mathbf{B}$ 
with sparser values than FastHALS-RP. 
Figure \ref{fig5.4b} shows that for both FastHALS and FastHALS-RP, the 
sparsity level of $\mathbf{B}$ decreases by increasing the value of the 
parameter $\beta$. In particular, for dense dataset the parameter $\beta$
does not introduce any relative variation between the two considered 
sparsity profiles. 

\subsection{20 Newsgroups dataset}
The 20 Newsgroups dataset \cite{lang1995newsweeder} is a collection of 
approximately 20000 text documents, evenly partitioned over 20 different 
newsgroups classes.
The analysis presented in this section is based on the application of the NMF 
methods on the frequency matrix associated with the 1000 most frequent words 
within the first 5000 samples. The resulting frequency matrix is sparse which
provides a completely different numerical framework compared to the one 
provided by Olivetti faces dataset. 
Following the same structure as in section \ref{sec5.1}, we first focus on the
reconstruction precision computed by 150 iterations of the NMF methods, with an 
arbitrary number of components $k = 60$ and with $r + r_{ov} = 72$ and $w = 9$ 
for the methods MU-RP, HALS-RP and FastHALS-RP.   

\begin{figure}[htbp]
  \includegraphics[width=0.95\textwidth]{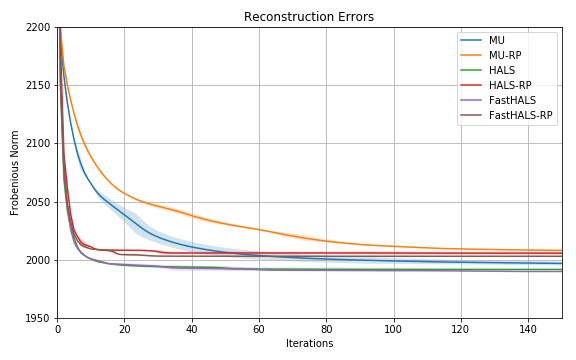}
  \caption{Evolution of the reconstruction error over iterations for the MU,
  MU-RP, HALS, HALS-RP, FastHALS, FastHALS-RP NMF methods when applied to the
  20 Newsgroup dataset.}
  \label{fig5.5}
\end{figure}

Figure \ref{fig5.5} shows that, also in the case of the 20 Newsgroup dataset, 
both MU and MU-RP are associated with a smaller convergence rate than the 
other NMF methods. Moreover, according to their iterative evolution profiles, 
these methods appear to be numerically more stable compared to the results
presented in Figure \ref{fig5.2} for the Olivetti faces dataset.
Moreover, Figure \ref{fig5.5} shows that NMF methods with random projections 
provide similar results to their standard counterparts in terms of data 
approximation error, with a few percent maximum deviation from the best 
performing algorithm FastHALS.

\begin{table}[htbp]
\centering
   \begin{tabular}{ |l||c|c||c| }
     \hline
     Method         & Complexity      & Time per update & Memory           \\
                  & [$10^6$ flops]  & [seconds]       & [$10^5$ float] \\
     \hline
     MU           & $2400$          & $0.068$         & $53$           \\      
     \hline
     MU-RP        & $150$           & $0.048$         & $12$           \\
     \hline
     HALS         & $2400$          & $2.678$         & $53$             \\
     \hline
     HALS-RP      & $430000$        & $4.775$         & $12$           \\
     \hline
     FastHALS     & $1300$          & $0.069$         & $53$             \\
     \hline 
     FastHALS-RP  & $98$            & $0.047$         & $12$           \\ 
     \hline
  \end{tabular}
\caption{Numerical complexity, average time per update and memory consumption 
         for the NMF of the 20 Newsgroups dataset.}
\label{table5.3}
\end{table}

From table \ref{table5.3} we can see that the theoretical value of the 
numerical complexity for the FastHALS-RP iterations is one order of magnitude
lower than the one of its unprojected counterpart FastHALS. This is reflected
by an averaged time per update $1.5$ times lower. Moreover, memory consumption
levels are reduced by a factor $\sim 3.5$ in this application case when using
random projections.

\begin{figure}[htbp]
\includegraphics[width=0.95\textwidth]{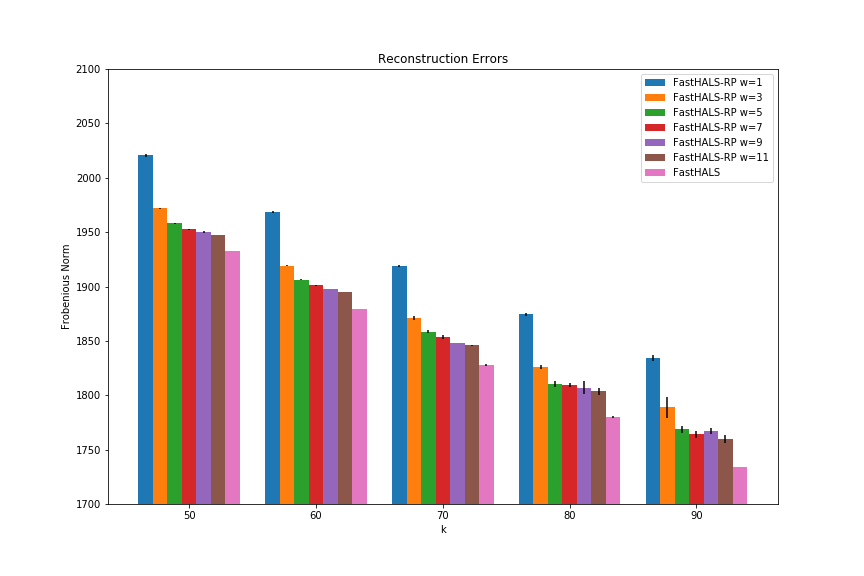}
\caption{Comparison of the median reconstruction errors resulting from the
    FastHALS and FastHALS-RP methods when considering $k = [50,60,70,80,90]$
    and $w = [1,3,5,7,9,11]$.}
        \label{fig5.6}
\end{figure}

In analogy with \ref{sec5.1}, we present a comparison between the results 
provided by FastHALS-RP and FastHALS to assess how the reconstruction errors 
are affected by different choices of the number of components $k$ and random 
projection parameter $w$. For this analysis we are considering
$k=[50,60,70,80,90]$ and $w=[1,3,5,7,9,11]$.
Figure \ref{fig5.6} shows that for this application the optimal choice is 
$w=11$, independently of the the number of components $k$. 
The discrepancy between this value for $w$ and the one found for the Olivetti 
faces dataset is due to the different decay rates $p$ of the singular value 
spectra $\sigma_j(X)^{-p}$ of the two dataset. In particular, for the 20 
Newsgroup dataset $p_{20\_Newsgroups} \sim 0.0012$ while for the Olivetti faces 
dataset $p_{Olivetti} \sim 0.0053$. 
With the ratio $\sim 20\%$ between $p_{20\_Newsgroups}$ and $p_{Olivetti}$, 
for the 20 Newsgroups dataset a higher value of $w$ is required to 
\emph{regularize} the distortion effects due to small singular values for the 
definition of the random projection operators $L$ and $R$.

\begin{figure}[htbp]
  \centering
  \includegraphics[width=0.95\textwidth]{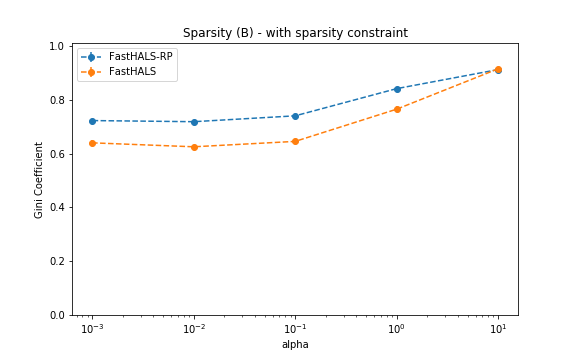}
  \caption{Gini coefficients of the $\mathbf{B}$ matrix computed by FastHALS 
           and FastHALS-RP with sparsity parameter 
           $\alpha=[0.001,0.01,0.1,1,10]$, number of components $k=60$ and
           random projection parameter $w = 9$.}
  \label{fig5.7a}
\end{figure}

\begin{figure}[htbp]
  \centering
  \includegraphics[width=0.95\textwidth]{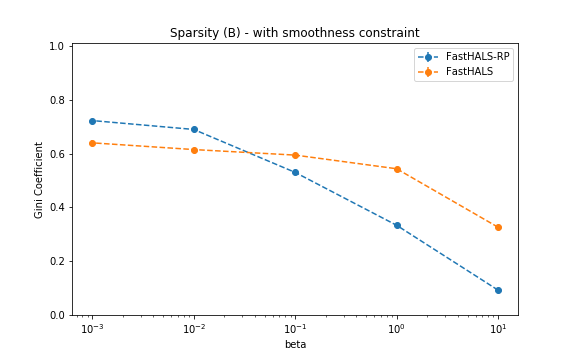}
  \caption{Gini coefficients of the $\mathbf{B}$ matrix computed by FastHALS 
           and FastHALS-RP with smoothness parameter 
           $\beta=[0.001,0.01,0.1,1,10]$, number of components $k=60$ and 
           random projection parameter $w = 9$.}
  \label{fig5.7b}
\end{figure}

Finally, we assess how the sparsity level of the factor matrix $\mathbf{B}$ 
changes for different choices of the parameters $\alpha$ and $\beta$ when 
keeping $k=60$ and $w=9$.

From Figure \ref{fig5.7a} and \ref{fig5.7b}, the sparsity and smoothness 
properties of $\mathbf{B}$ show coherent behaviors to the one presented 
for the Olivetti Faces dataset. In particular, the Gini coefficients 
profiles for both the solutions provided by FastHALS and FastHALS-RP grow 
as function of $\alpha$ and decrease for larger values of $\beta$. 
 
Figure \ref{fig5.7a} shows that, for the same set of values for $\alpha$, 
FastHALS-RP provides solutions which are generally sparser than FastHALS. 
Moreover, the two methods seem to react differently to different choice of 
the parameter $\beta$, as shown in Figure \ref{fig5.7b}. In particular, over 
the same range of $\beta$ values we assessed that FastHALS-RP provides a 
stronger variations in the profile of Gini coefficients than FastHALS. 

In conclusion, both Figure \ref{fig5.7a} and \ref{fig5.7b} shows that for a 
sparse dataset and without imposing any sparsity and smoothness constraint, 
the solutions for $\mathbf{B}$ provided by FastHALS-RP are generally sparser 
than the FastHALS ones. All those properties are a direct consequence of the 
application of Random Projections on a sparse dataset. 

\section{Conclusion}

This paper was aimed at proposing FastHALS-RP, a novel NMF method capable 
of combining the fastest existing NMF algorithm FastHALS with a data 
dimensionality reduction scheme based on Random Projections.
Our new algorithm outperforms the state of the art NMF methods both in 
computational efficiency and memory consumption. \\
With an appropriate choice of the parameters $k$ and $w$, the desired number of
components and the random projection power iteration parameter respectively, it
has been shown  that FastHALS-RP provides results with a data approximation
precision level very close to the standard FastHALS method. In particular, we
showed that the optimal choice for the parameter $w$ strongly depends on the
decay rate of the singular value spectrum of the data.\\
We showed that sparsity and smoothness constraints on the factor matrix
$\mathbf{B}$ can easily be introduced into the FastHALS-RP optimization
problem. The conclusion of our analysis focused on assessing how the sparsity
properties of the factor matrix $\mathbf{B}$, computed by FastHALS and
FastHALS-RP, are influenced by different choices of $\alpha$ and $\beta$,
the sparsity and smoothness coefficients in the likelihood function respectively.
We found that for dense datasets and in absence of any sparsity and smoothness 
constraints for $\mathbf{B}$, FastHALS and FastHALS-RP provide solutions with 
similar sparsity property, while for sparse dataset, FastHALS-RP provides 
sparser solutions. For this reason, the effect of the two parameters $\alpha$ 
and $\beta$ is different in the two scenarios and needs to be properly 
evaluated from case to case.\\
Some theoretical aspects regarding the optimal choice of the parameters 
$w$, $r + r_{ov}$, $\alpha$ and $\beta$ still need to be investigated as 
well as the global convergence property of the algorithm.\\
While we could substantially reduce the memory footprint and speed up
computation compared to existing methods, our algorithm still relies on being
computed on a single node. A next step towards large-scale applications on
distributed systems could now be to extend our algorithm with the map-reduce
scheme proposed by \cite{yu2014parallel}. We expect that such an approach
will profit twice from our algorithm since we will not only see a speed-up on
all nodes individually, but, in addition, the data transfer between
nodes will be reduced which will lead to an additional speed-up.\\

We implemented the proposed method in Python and the source code will soon be publicly available.

\section*{Acknowledgements}
Our work was supported by the two SNF Synergia Grant (\emph{EUCLID: 
high-precision cosmology in the dark sector}). We Would like to thank 
Martin Melchior, André Csillaghy and Roman Bolzern for helping and 
providing comments that greatly improved the quality of this manuscript.

\bibliography{biblio}{}
\bibliographystyle{plain}

\end{document}